\documentclass[11pt,a4paper]{report}
\usepackage[utf8]{inputenc}
\usepackage[english]{babel}
\usepackage{amsmath}
\usepackage{amsfonts}
\usepackage{amssymb}
\usepackage{graphicx}

\newtheorem{definition}{Definition}
\newtheorem{open problem}{Open Problem}

\author{Andrés García and Juan Andrés Roteta Lannes}%

\begin{document}

\title{A generalized initial value problem for ODE's}
\maketitle

\abstract{
This paper introduces, up to the author's knowledge, for the first time the generalized initial value problem. In this problem, given an ordinary differential equation defined in some set, the initial conditions are mapped to a subset of that domain.
This generalization allows some systems to radically change their properties, in fact some obstructions to asymptotic stability are get rid off by mapping the initial conditions set using modularity.
}

\section{Introduction}

The initial value problem for ordinary differential equations (ODE's) is a basic and the primal problem to consider in any dynamic or control system (see for instance \cite{Chicone2006}). Among the different possible trajectories that an ODE can take, asymptotic stability is always a mater of research interest among the control community (see for instance \cite{Bacciotti2005}).

To date, even with many available techniques, stability, asymptotic stability or even controller design for non-linear systems is not an easy task (see for instance \cite{Sastry1999}). The case of asymptotic stability offers a very challenging problem, even using the well studied Lyapunov's method providing only a sufficient condition leaving the user the construction of such a function. 

In the case of control systems (systems modelled with ODE's but with the inputs as free parameters \cite{Sastry1999}), the problem becomes involved since existence theorems are needed in order to determine beforehand the possibility of asymptotic stabilization (see for instance \cite{Brockett1983} and \cite{Garcia2018}).

The paper in \cite{Garcia2020} a smooth controller to render the origin of an unicycle robot asymptotically stable was derived. This clearly seems to contradict the Brockett's condition, which for this kinematic system proved the impossibility of such a controller.

However, a closer inspection of the controller designed reveals that the modularity in one state-space variable: $\theta(t)$ was invoked to mapp the set of initial conditions. In this case, only the initial conditions are mapped from the complete $\Re^{3}$ set for the value $\theta(0)=0$.

In this paper, motivated by the results in \cite{Garcia2020}, a generalized initial value problem is defined in order to present a new perspective but also an open problem to define theorems and tools to determine the possibility for asymptotic stabilizing controllers.

This paper is organized as follows: In Section \ref{GIVP} the mathematical set-up and the definition of the generalized initial value problem is presented along with an open problem, whereas Section \ref{Conclusions} presents some conclusions and future work.

\section{The initial conditions mapping: different trajectories}\label{GIVP}

\begin{figure}
		\centering
		\includegraphics[width=0.5\textwidth]{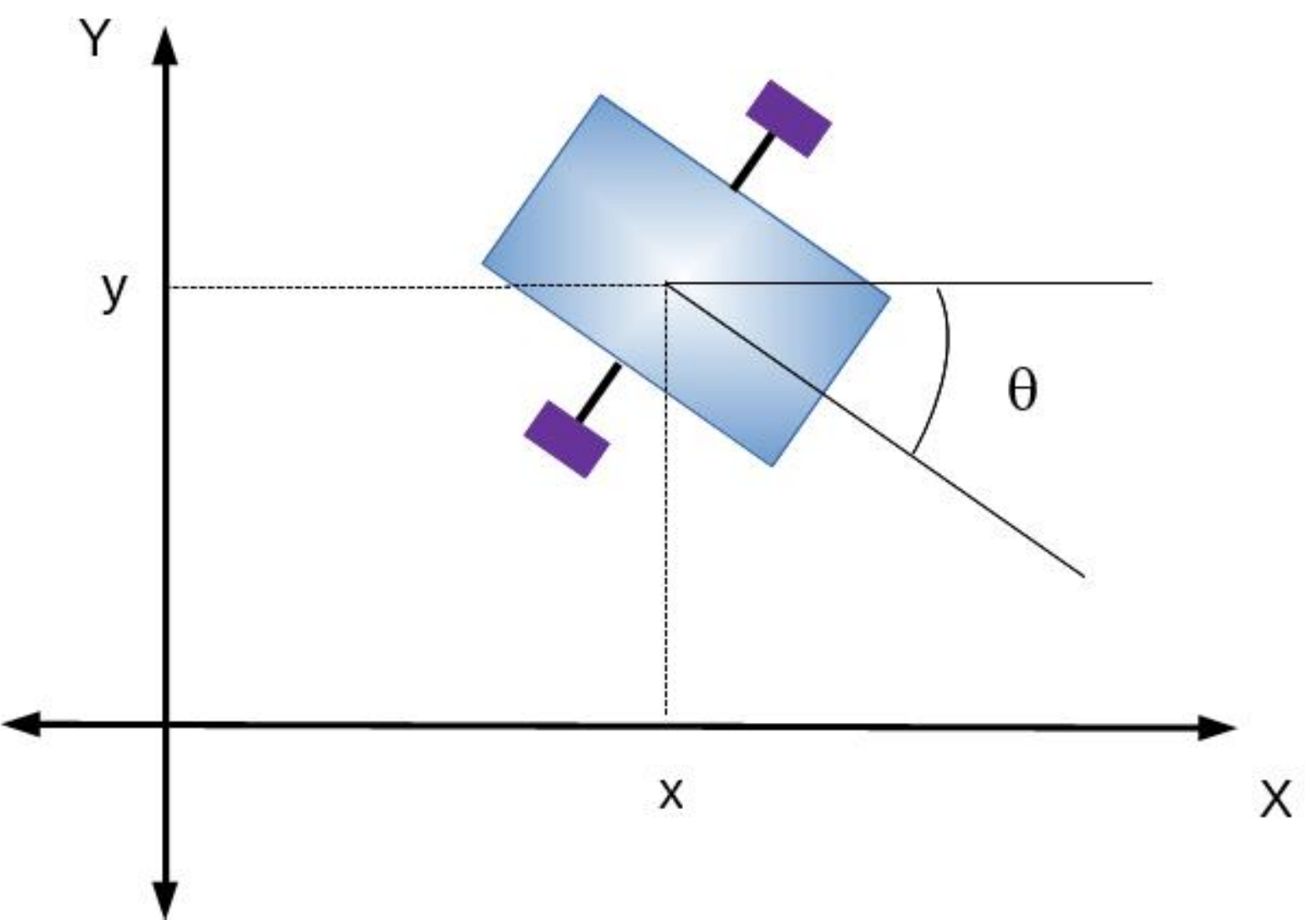}
		\caption{Unicycle-like robot kinematic's model} \label{Robot-Model}
	\end{figure}

The problem studied in this paper considers a generalization of the well-known initial value problem (see for instance \cite{Chicone2006}):

\begin{definition}[Generalized initial value problem]\label{GIVP: Definition}
Given an ODE $\dot{x}=f(x),\quad x\in\Omega\subset\Re^{n}$ and an initial condition $x(0)=\phi(x_{0}),\quad x_{0}\in\Omega$ along with $\phi:\Omega\rightarrow \Omega^{*}\subset\Omega$.

Finding trajectories (see \cite{Chicone2006} for details on trajectories definitions) satisfying this definition is called \textit{generalized initial value problem}. 
\end{definition}

where $\dot{x}$ means the time derivative. A motivation for such a generalization comes from a well-known kinematic model for mobile robots (see for instance \cite{Garcia2012} and Figure \ref{Robot-Model}):

	\begin{equation}\label{Kinematic model}
		\begin{bmatrix}
			\dot{x}\\
			\dot{y}\\
			\dot{\theta}
		\end{bmatrix}=f(x,u)=
		\begin{bmatrix}
			cos(\theta) && 0\\
			sin(\theta) && 0\\
			0 && 1
		\end{bmatrix} \cdot
		\begin{bmatrix}
			u_{1}\\
			u_{2}
		\end{bmatrix}
	\end{equation}

where $\{u_{1},u_{2}\}$ are the control inputs or control variables. It is well known that asymptotic stability (see for instance ) can not be realized by means of smooth controllers (see for instance \cite{Brockett1983} and \cite{Garcia2018}).

For the sake of completeness the conditions to drive this model to the origin using smooth controllers can be readily written as follows (see\cite{Garcia2018}):

	\begin{equation*}
		f(V_{1},u(V_{1}))=\alpha \cdot f(V_{2},u(V_{2})), \quad \alpha \neq 0,k\quad \left\{V_{1},V_{2}\right\} \in 
		\delta_{0}		
	\end{equation*}

where $\delta_{0}$ is a neighbour of the origin. This condition is equivalent to look for constant direction regions:

	\begin{equation*}
		\frac{f(x,u(x))}{\Vert f(x,u(x))\Vert}=constant,\quad \forall x\in\
		\delta_{0}
	\end{equation*}

Checking this condition in equation (\ref{Kinematic model}):
	
	\begin{equation*}
	 \begin{cases}
		\frac{cos(\theta)}{sin(\theta)}=\rho_{1}\\
		\frac{cos(\theta)\cdot u_{1}}{u_{2}}=\rho_{2}\\
		\frac{sin(\theta)\cdot u_{1}}{u_{2}}=\rho_{3}
	  \end{cases}
	\end{equation*}

where $\{\rho_{1},\rho_{2},\rho_{3}\}$ are constants. Clearly, the condition $\frac{cos(\theta)}{sin(\theta)}=\rho_{1}$ means a straight line excluding all the rest as asymptotically stable trajectories.

However, the controller defined in \cite{Garcia2020} renders the model (\ref{Kinematic model}) globally asymptotically stable:
	
	\begin{equation}\label{Controller: Bessel}
		\begin{bmatrix}
			u_{1}\\
			u_{2}
		\end{bmatrix}=
		\begin{bmatrix}
			\sum_{i=1}^{N} (2 cdot i+1) \cdot a \cdot C_{i} \cdot J_{i}(\theta)
			\cdot \theta^{i+1}\\
			a \cdot \theta, \quad a<0
		\end{bmatrix}
	\end{equation}		

where $N \in \mathbf{N}$ is an arbitrary number and $\{C_{i},a<0\}$ some constants. Moreover, an important initial conditions mapping is needed:

	\begin{equation}\label{Initial conditions}
		\theta(0)=
		\begin{bmatrix}
			\theta_{0}\neq 0,\quad \theta_{0}\in \Re\\
		    2 \cdot \pi,\quad \theta_{0}=0
		\end{bmatrix}
	\end{equation}

Notice that without this initial conditions mapping, the system (\ref{Kinematic model}) along with controller (\ref{Controller: Bessel}) is not asymptotically stable:

	\begin{equation*}
		\begin{bmatrix}
			\dot{x}\\
			\dot{y}\\
			\dot{theta}
		\end{bmatrix}=
		\begin{bmatrix}
			cos(\theta) && 0\\
			sin(\theta) && 0\\
			0 && 1
		\end{bmatrix} \cdot
		\begin{bmatrix}
			\sum_{i=1}^{N} (2 \cdot i+1) \cdot a \cdot C_{i} \cdot J_{i}(\theta)
			\cdot \theta^{i+1}\\
			a \cdot \theta
		\end{bmatrix}
	\end{equation*}

Considering the initial condition $\theta(0)=0$:

	\begin{equation*}
		 \dot{\theta}=a \cdot \theta \Leftrightarrow \theta(t)=0 \quad \forall t
		 \in\Re^{+}
	\end{equation*}

Clearly, this is not asymptotically stable in $\{x(t),y(t)\}$ to the origin:

	  \begin{eqnarray*}	  
		\begin{bmatrix}
			\dot{x}\\
			\dot{y}
		\end{bmatrix}=
		\begin{bmatrix}
			cos(\theta)\\
			sin(\theta)
		\end{bmatrix} \cdot
		\sum_{i=1}^{N} (2 \cdot i+1) \cdot a \cdot C_{i} \cdot J_{i}(\theta)
		\cdot \theta^{i+1}\\
		\theta(t)=0		
	  \end{eqnarray*}
	
That is:

	\begin{equation*}
		\begin{bmatrix}
			x-x(0)\\
			y-y(0)
		\end{bmatrix}=
		\begin{bmatrix}
			1\\
			0
		\end{bmatrix} \cdot
		\sum_{i=1}^{N} (2 \cdot i+1) \cdot a \cdot C_{i} \cdot J_{i}(0)
		\cdot 0^{i+1}
	\end{equation*}
	
Finally:

	\begin{equation*}
		\begin{bmatrix}
			x-x(0)\\
			y-y(0)
		\end{bmatrix}=
		\begin{bmatrix}
			0\\
			0
		\end{bmatrix} 
	\end{equation*}

This conclusions showed that mapping the initial condition set (not the entire trajectories domain) changes completely the behaviour of the orbits.

\begin{figure}
		\centering
		\includegraphics[width=0.4\textwidth]{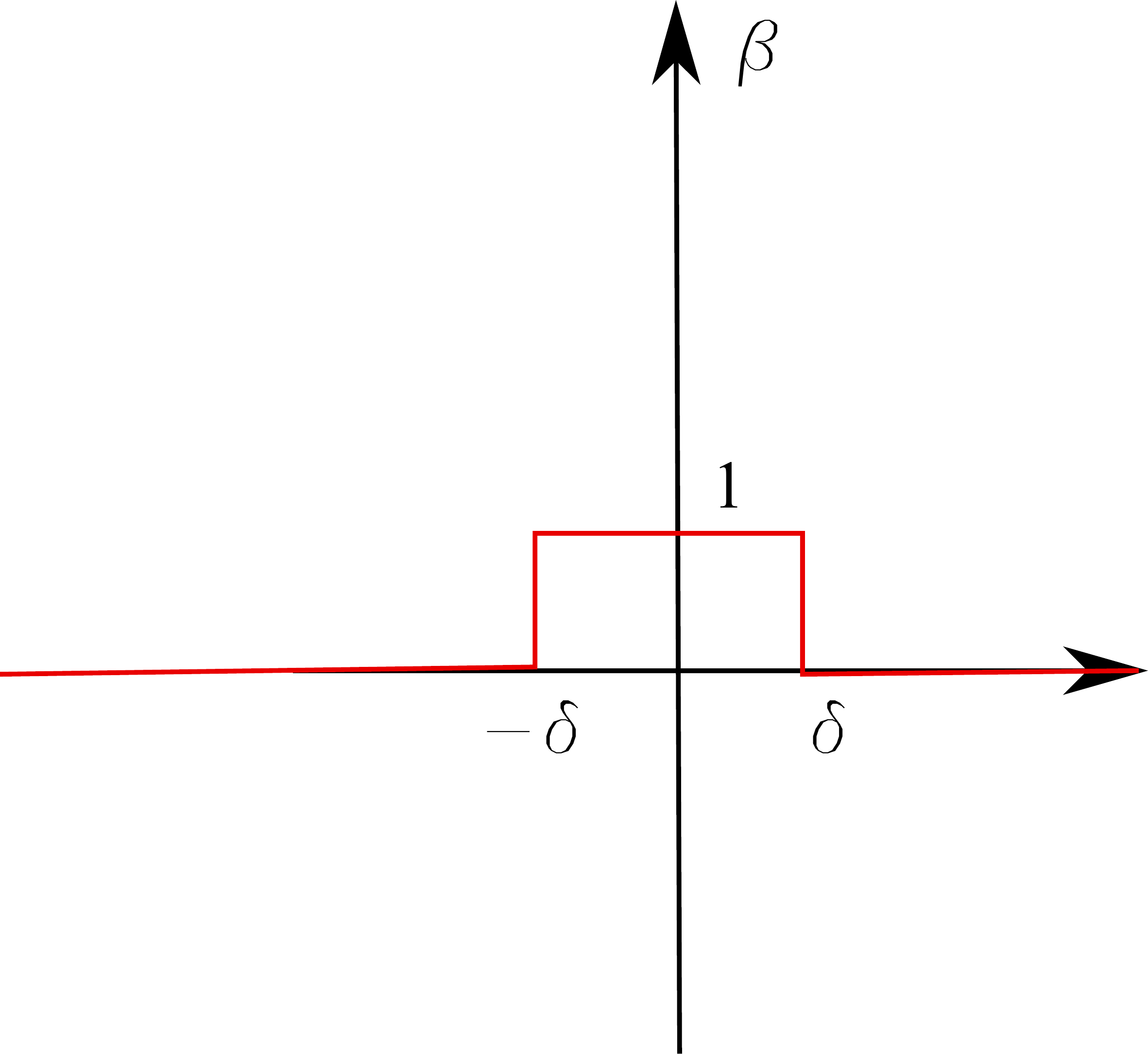}
		\caption{Function $\beta(\theta_{0})$} \label{beta}
	\end{figure}

It is worth to notice that the general mapping definition (\ref{GIVP: Definition}) for this case can be written as follows:

	\begin{equation*}
		\theta(0)=
		\begin{bmatrix}
			\theta_{0}\neq 0,\quad \theta_{0}\in \Re\\
		    2 \cdot \pi,\quad \theta_{0}=0
		\end{bmatrix} \Leftrightarrow \theta(0)=\theta_{0}+2\cdot \pi 
		\beta(\theta_{0})
	\end{equation*}

where the function $\beta(\theta_{0})$ is depicted in Figure \ref{beta}, with $\beta \longrightarrow 0$.

It turns out that the classical results, theorems and method for the initial value problem with $\phi$ the identity function can not be applied to the generalized definition \ref{GIVP: Definition}:

\begin{open problem}
 Given a control system: $\dot{x}=f(x,u), x\in\Re^{n}$ with control inputs: $u\in\Re^{m}, m \leq n$ derive conditions for asymptotic stability in the case of generalized initial value problem.
\end{open problem}

\textbf{Note:} For this generalized problems, it is necessary to identify the subset of initial conditions  that can be mapped onto an equivalent values of another ODE's definition domain subset. In other words, the existence's possibility and its physical meaning for the mapping function $\phi$ in (\ref{GIVP: Definition}) must be established beforehand.

\section{Conclusions}\label{Conclusions}

In this paper and up to the authors knowledge, a new formulation to the classical initial value problem for ODE's is presented as a generalized initial value problem.

The motivation for this generalization comes from the smooth controller found in \cite{Garcia2020} in the case of a unicycle robot model that clearly would violate the well-known Brockett's condition.

The key to this result lies into the initial value mapping that changes completely the trajectories' behaviour allowing asymptotic stability even using smooth controllers.

This new result and formulation open the road to reconsider smooth controller's in the cases where Brockett's condition forbid such designs but, at the same time, open the necessity for new theorems and methods to establish the existence of smooth controllers for asymptotic stability in the case of generalised initial value problem.

As a future work, more complex models using ODE's will be analysed using these ideas in order to provide more examples were initial value mapping is useful providing asymptotic stability.

\section{Acknowledgments}
The authors would like to acknowledge Universidad Tecnol\'{o}gica Nacional, Facultad Regional Bah\'{i}a Blanca. and Comisión de Investigaciones Coentíficas (CIC) and Universidad Tecnológica Nacional-FRBB.

\end{document}